\documentclass[12pt]{article}
\usepackage{amssymb}
\usepackage{color}
\usepackage{cite}

\begin{document}
\newcommand{\lap}{\mbox{$\bigtriangleup$}}
\newcommand{\grad}{\mbox{$\bigtriangledown$}}
\newcommand{\be}{\begin{equation}}
\newcommand{\ee}{\end{equation}}
\newcommand{\kernel}{\frac{1}{ |x-y|^{n-\alpha} } }
\newcommand{\ind}{ (n+\alpha)/(n-\alpha) }
\newcommand{\indi}{ \frac{n+\alpha}{n-\alpha} }
\newcommand{\bound}{ (1+|x|)^{\alpha-n} }
\newcommand{\bounds}{\frac{1}{ (1+|x|)^{n-\alpha} } }
\newcommand{\solution}{ (\frac{t}{t^{2}+|x-x_{0}|^{2} } )^{ (n-\alpha)/2 } }
\newcommand{\sol}{ u_{t,x_{0}}(x) }
\newcommand{\half}{ \frac{1}{2} }
\newcommand{\hsp}{ \Sigma_{\lambda} }
\newcommand{\hspn}{ \Sigma_{\lambda}^{-} }
\newcommand{\hspp}{ \Sigma_{\lambda}^{+} }
\newcommand{\dom}{D(\lambda, \delta, R_{1})}
\newtheorem{mthm}{Theorem}
\newtheorem{mcor}{Corollary}
\newtheorem{mpro}{Proposition}
\newtheorem{mfig}{figure}
\newtheorem{mlem}{Lemma}
\newtheorem{mdef}{Definition}
\newtheorem{mrem}{Remark}
\newtheorem{mcon}{Conjecture}
\newtheorem{mpic}{Picture}
\newtheorem{rem}{Remark}[section]
\newcommand{\ra}{{\mbox{$\rightarrow$}}}
\newtheorem{thm}{Theorem}[section]
\newtheorem{pro}{Proposition}[section]
\newtheorem{lem}{Lemma}[section]
\newtheorem{defi}{Definition}[section]
\newtheorem{cor}{Corollary}[section]

\title{\bf  Existence, Uniqueness of Positive Solution to a Fractional Laplacians with Singular Nonlinearity}

\author{
Yanqin Fang \thanks{The author is supported by National Natural Sciences Foundations of China, grant No. 11301166, and Young Teachers Program of Hunan University 602001003}
}
\date{\today}
\maketitle

\begin{abstract}
In this paper we prove the existence and uniqueness of positive classical solution of the fractional Laplacian with singular nonlinearity
in a smooth bounded domain with zero Drichlet boundary conditions.
By the method of sub-supersolution, we derive the existence of positive classical solution to the
approximation problems.  In order to obtain
the regularity, we first establish the existence of weak solution for the fraction Laplacian. Thanks to \cite{XY}, the regularity follows from the boundedness of
weak solution.

\end{abstract}

\bigskip


\section{Introduction}
We are concerned in this paper the existence and uniqueness of classical solution to
the nonlinear boundary value problems involving the fractional Laplacian
\begin{equation}\label{01}
\left\{
  \begin{array}{ll}
(-\Delta)^{s} u =u^{-p},&x\in\Omega,\\
u(x)>0,&x\in\Omega,\\
u(x)=0,&x\in\partial\Omega,
  \end{array}
\right.
\end{equation}
where $\Omega$ is a smooth bounded domain of
$\mathbb{R}^{N}$, $0<s<1$, $0<p<1$. The operator $(-\Delta)^{s}$ is defined by the eigenvalues and eigenfunctions of
the original operator $-\Delta$. It can be understood as the nonlocal
version of the Brezis-Nirenberg problem \cite{BN}. We are looking for classical solutions of (\ref{01}), i.e., solutions which at least
belong to the class $C_{loc}^{2,\alpha}(\Omega)$, $0<\alpha<1$.

In recent years, there has been tremendous interested in developing the problems related to Laplacian:
\begin{equation}\label{14030709}
-\Delta u=u^{-p}+f(x,u),\;\;p>0,\;\;\mbox{in}\;\Omega.
\end{equation}
Let us mention the pioneering papers of Crandall, Rabinowitz and Tartar \cite{CRT} and Stuart \cite{S}. For example, in \cite{CRT}, the authors considered nonlinear elliptic boundary value problems of the form
\begin{equation}\label{14030901}
\left\{
  \begin{array}{ll}
   -\sum_{i,j=1}a_{ij}(x)\frac{\partial^{2}u}{\partial x_{i}\partial x_{j}}+\sum_{i=1}^{n}b_{i}(x)\frac{\partial u}{\partial x_{i}}+c(x)u=g(x,u),&x\in\Omega\\
u(x)=0,&x\in\partial\Omega.
  \end{array}
\right.
\end{equation}
They presented the following proposition.
\begin{mpro}\cite{CRT}
If $g$ satisfies

$(g_{1})$
$$
\lim_{r\rightarrow 0^{+}}g(x,r)=\infty \;\mbox{uniformly}\;\mbox{for}\; x\in\bar{\Omega},
$$
and

$(g_{2})$
$$
g(x,r)\;\mbox{is}\;\mbox{nonincreasing}\;\mbox{in}\;r\in (0,\infty)\;\mbox{for}\;x\in\bar{\Omega},
$$
then (\ref{14030901}) possesees a unique classical solution $u\in C^{2}(\Omega)\cap C(\bar{\Omega})$ with $u>0$ in $\Omega$.
\end{mpro}
Their proofs are based on
finding sub- and supersolutions for approximate problems together with Sobole embedding and appropriate a priori estimates.
In this work we are inspired to follow their ideas---employing sub- and supersolutions and some new technics for approximate problems related to fractional Laplacian.

In \cite{CD}, Canino and Degiovanni provided a variational approach to singular semilinear elliptic equations of the form (\ref{14030709}).
In a recent paper \cite{CGB}, Canino, Grandinetti and Sciunzi considered positive solutions to the singular semilinear elliptic equation. By the moving plane method, they
deduced symmetry and monotonicity properties of solutions. For more results related to (\ref{14030709}), we refer to \cite{Bo,BO,Ca,GNN,HSS,LS,LM}

The fractions of the Laplacian are the infinitesimal generators of L\'{e}vy stable diffusion processes \cite{A}.
It appears in diverse areas including physics, biological modeling
and mathematical finances. There have been extensive study the partial differential equations involving the fractional
Laplacian. Recently, after the work of  Caffarelli and Silvestre \cite{CS}, several authors have studied an
equivalent definition of the operator $(-\Delta)$ in a bounded domain with zero Dirichlet boundary. Results relating to these problems can be found in \cite{BCP,CaSi,CT,CDDS,SV,Si,ST,TZ }.

Motivated by some results found in \cite{CaSi, CT,CD,CDDS, CRT, Si, XY}, a natural question arises whether the existence and uniqueness of classical solutions can be
obtained in fraction Laplacian with singular nonlinearity. As far as I known, little work has been done in fractional Laplacian
with such nonlinearity. The nonlocal property of fractional Laplacian and singularity make it difficult to handle.

Let $\{\lambda_{i},\varphi_{i}\}$ be the eigenvalues and corresponding eigenfunctions of $(-\Delta)$ in $\Omega$ with
zero Dirichlet boundary condition:
$$
\left\{
  \begin{array}{ll}
  -\Delta \varphi_{i}=\lambda_{i}\varphi_{i}&\mbox{in}\;\Omega,\\
\varphi_{i}=0 &\mbox{on}\;\partial\Omega,
  \end{array}
\right.
$$
It is well known that $0<\lambda_{1}<\lambda_{2}\leq \lambda_{3}\leq\cdots\leq\lambda_{k}\rightarrow+\infty$, and
the corresponding eigenfunctions are orthonormal, i.e,
$$
\int_{\Omega}\varphi_{i}\varphi_{j}dx=\delta_{ij}=\left\{
                                                    \begin{array}{ll}
                                                     1,& i=j,\\
0,&i\neq j.
                                                    \end{array}
                                                  \right.
$$
The operator $(-\Delta)^{s}$ is defined for any $u\in C^{\infty}_{c}(\Omega)$ by
$$
(-\Delta)^{s}u=\sum_{i=1}^{\infty}\lambda^{s}_{i}u_{i}\varphi_{i},
$$
where
$$
u=\sum_{i=1}^{\infty}u_{i}\varphi_{i},\;\;\mbox{and}\;u_{i}=\int_{\Omega}u\varphi_{i}dx.
$$
The operator can be extended by density for $u$ in the Hilbert space
$$
H^{s}_{0}(\Omega)=\left\{u\in L^{2}(\Omega):\;\|u\|^{2}_{H^{s}_{0}(\Omega)}=\sum_{i=1}^{\infty}\lambda^{s}_{i}|u_{i}|^{2}<+\infty\right\}.
$$
Note that then
$$
\|u\|_{H^{s}_{0}(\Omega)}=\|(-\Delta)^{\frac{s}{2}}u\|_{L^{2}(\Omega)}.
$$

Now, we give the definition of weak solution of (\ref{01}).
\begin{mdef}
We say that $u\in H^{s}_{0}(\Omega)$ is a solution of (\ref{01}) if the identity
\begin{equation}
\int_{\Omega}(-\Delta)^{\frac{s}{2}}u(-\Delta)^{\frac{s}{2}}\varphi dx=\int_{\Omega}u^{-p}\varphi dx
\end{equation}
holds for every function $\varphi\in H^{s}_{0}(\Omega)$.
\end{mdef}

As we all know that the operator $(-\Delta)^{s}$ can be realized as the
boundary operator of a suitable extension in the
half-cylinder $\Omega\times(0,\infty)$. Such
an interpretation was demonstrated in \cite{CS} for the fractional Laplacian in $\mathbb{R}^{N}$.

Indeed, let us define
\begin{equation}\label{14011701}
\mathcal{C}=\Omega\times (0,+\infty), \;\;\partial_{L}\mathcal{C}=\partial\Omega\times[0,+\infty).
\end{equation}
We write points in the cylinder $\mathcal{C}$ by $(x,y)\in\mathcal{C}=\Omega\times (0,+\infty)$.
Given $s\in (0,1)$, consider the space $H^{1}_{0,L}(y^{1-2s})$ of measurable functions $U:\mathcal{C}\rightarrow \mathbb{R}$
such that $U\in H^{1}(\Omega\times(s,t))$ for all $0<s<t<+\infty$, $U=0$ on $\partial_{L}\mathcal{C}$
and for which the following norm is finite:
\begin{equation}\label{14011203}
\|U\|^{2}_{H^{1}_{0,L}(y^{1-2s})}=\int_{\mathcal{C}}y^{1-2s}|\nabla U|^{2}dxdy.
\end{equation}

Recall that if $\Omega$ is a smooth bounded domain, it is verified that
$$
H^{s}_{0}(\Omega)=\{\;u=\mbox{tr}|_{\Omega\times\{0\}}U:\;U\in H^{1}_{0,L}(\mathcal{C})\}.
$$
in \cite{CS,Ta}.
In order to study the existence and uniqueness of positive solution of (\ref{01}), we may consider
the fractional harmonic extension of a function $u$ defined in $\Omega$: (See \cite{BCP,CT,CS,ST,Ta})
\begin{equation}\label{14011202}
\left\{
  \begin{array}{ll}
   \mbox{div}(y^{1-2s}\nabla U)=0,&\mbox{in}\;\mathcal{C},\\
   U=0,&\mbox{on}\;\partial_{L}\mathcal{C},\\
   y^{1-2s}\frac{\partial U}{\partial \nu}=u^{-p},&\mbox{in}\;\Omega\times\{0\},
  \end{array}
\right.
\end{equation}
where $\nu$ is the unit outer normal to $\Omega\times\{0\}$. If $U\in H^{1}_{0,L}(y^{1-2s})$ satisfies (\ref{14011202}),
then the trace $U$ on $\Omega\times\{0\}$ of the function $U$ will be a solution of problem (\ref{01}) in the weak sense.

The main result is the following theorem.

\begin{mthm}\label{mthm14011801}
Assume $\Omega\subset\mathbb{R}^{N}$ is a smooth bounded domain, $0<s<1$, $0<p<1$. Then for some $0<\alpha<1$, there exist a unique classical solution  $u\in C_{loc}^{2,\alpha}(\Omega)$ to
\begin{equation}\label{0310}
\left\{
  \begin{array}{ll}
(-\Delta)^{s} u =u^{-p},&x\in\Omega,\\
u(x)>0,&x\in\Omega,\\
u(x)=0,&x\in\partial\Omega.
  \end{array}
\right.
\end{equation}
\end{mthm}

The proof of Theorem \ref{mthm14011801} is given in section 2.

\section{Proof of Theorem \ref{mthm14011801}}
Before proof of Theorem, we state some lemmas.
\begin{lem}\cite{CDDS}\label{mlem14022101}
Let $h\in (H^{s}_{0}(\Omega))^{\ast}$. Then, there is a unique solution to the problem:
$$
\mbox{find}\;u\in H^{s}_{0}(\Omega)\;\mbox{such that} \;(-\Delta)^{s}u=h.
$$
Moreover $u$ is the trace of $U\in H^{1}_{0,L}(y^{1-2s})$, where $v$ is the unique solution to
\begin{equation}\label{14011805}
\left\{
  \begin{array}{ll}
    \mbox{div}(y^{1-2s}\nabla U)=0&\mbox{in}\;\mathcal{C}\\
U=0&\mbox{on}\;\partial_{L}\mathcal{C}\\
-\lim_{y\rightarrow 0}(y^{1-2s}\frac{\partial U}{\partial y})=c_{N,s}h&\mbox{on}\;\Omega\times\{0\}
  \end{array}
\right.
\end{equation}
where $c_{n,s}>0$ is a constant depending on $N$ and $s$ only, equation (\ref{14011805}) is understood in the sense that
$U\in H^{1}_{0,L}(y^{1-2s})$ and
\begin{equation}
c_{N,s}\langle h, \mbox{tr}_{\Omega}(v)\rangle_{(H^{s}_{0}(\Omega))^{\ast},H^{s}_{0}(\Omega)}=\int_{\mathcal{C}}y^{1-2s}\nabla U\nabla v dxdy, \;\;\forall\;v\in H^{1}_{0,L}(y^{1-2s}),
\end{equation}
where $\langle h, \mbox{tr}_{\Omega}(v)\rangle_{(H^{s}_{0}(\Omega))^{\ast},H^{s}_{0}(\Omega)}$ is the duality pairing between
$(H^{s}_{0}(\Omega))^{\ast}$ and $H^{s}_{0}(\Omega)$. The constant $c_{n,s}$ is the same constant appearing in (\ref{14011805}).
\end{lem}
\begin{lem}\cite{CaSi,CDDS}\label{mlem14022102}
Let $h\in (H^{s}_{0}(\Omega))^{\ast}$ and $U\in H^{1}_{0,L}(y^{1-2s})$ denote the solution of (\ref{14011805}). Then, for any $\omega\subset\subset \Omega$, $R>0$, we have

$(i)$ If $h\in L^{\infty}(\Omega)$, then $U\in C^{\beta}(\omega\times [0,R])$, for any $\beta\in (0,\min(1,2s))$,

$(ii)$ If $h\in C^{\beta}(\Omega)$, then
\begin{enumerate}
  \item $U\in C^{\beta+2s}(\omega\times [0,R])$ if $\beta+2s<1$,
  \item $\frac{\partial U}{\partial x_{i}}\in C^{\beta+2s-1}(\omega\times [0,R])$ if $1<\beta+2s<2$, $i=1,\cdots,N$,
  \item  $\frac{\partial^{2} U}{\partial x_{i}\partial x_{j}}\in C^{\beta+2s-2}(\omega\times [0,R])$ if $2<\beta+2s$, $i,j=1,\cdots,N$,
\end{enumerate}
\end{lem}

\begin{pro}\cite{XY}\label{mpro14030701}
Suppose
$$
|f(t)|\leq C(1+|t|^{p})
$$
holds. Let $u\in H^{s}(\Omega)$ be a solution of
$$
\left\{
  \begin{array}{ll}
    (-\Delta)^{s}u=f(u)&\mbox{in}\;\Omega,\\
u=0&\mbox{on}\;\partial\Omega.
  \end{array}
\right.
$$
Then $u\in L^{\infty}(\Omega)$ if $1<p<\frac{n+2s}{n-2s}$; and $u\in L^{\infty}_{loc}(\Omega)$ if $p=\frac{n+2s}{n-2s}$.
Consequently, $u\in C^{2,\alpha}_{loc}(\Omega)$ for some $0<\alpha<1$.
\end{pro}

\textbf{Proof of Theorem \ref{mthm14011801}.} We first establish the existence of the approximate problems
\begin{equation}\label{02}
\left\{
  \begin{array}{ll}
  (-\Delta)^{s} u_{\epsilon}=(\epsilon+u_{\epsilon})^{-p}&\mbox{in}\;\Omega\\
u_{\epsilon}>0&\mbox{in}\;\Omega\\
u_{\epsilon}=0&\mbox{on}\;\partial\Omega
  \end{array}
\right.
\end{equation}
for $0<\epsilon<1$, and then show the convergence of $u_{\epsilon}$ as $\epsilon\rightarrow 0^{+}$ to a solution $u$.

For any fixed $\epsilon>0$, a function $v_{\epsilon}$ is called a subsolution
of (\ref{02}) if $v_{\epsilon}\in C_{loc}^{2,\alpha}(\Omega)$ and
$$
\left\{
  \begin{array}{ll}
   (-\Delta)^{s} v_{\epsilon}\leq(\epsilon+v_{\epsilon})^{-p}&\mbox{in}\;\Omega,\\
v_{\epsilon}=0&\mbox{on}\;\partial\Omega,
  \end{array}
\right.
$$
and $w_{\epsilon}$ is called a supersolution if $w_{\epsilon}\in C_{loc}^{2,\alpha}(\Omega)$ and
$$
\left\{
  \begin{array}{ll}
     (-\Delta)^{s} w_{\epsilon}\geq(\epsilon+w_{\epsilon})^{-p}&\mbox{in}\;\Omega,\\
w_{\epsilon}=0&\mbox{on}\;\partial\Omega.
  \end{array}
\right.
$$

In order to obtain the existence of solutions, we introduce the following lemma. It can be derived from
lemmas 4.2 and 4.3 in \cite{BCP}.
\begin{lem}\label{mlem14022110}
Assume there exist a  subsolution $v$ and a supersolution $w$ to
problem (\ref{02}). Then there exists a classical solution $u$ satisfying
$v\leq u\leq w$.
\end{lem}

Now we establish
\begin{lem}\label{14031105}
Let $0<\epsilon_{0}<1$. If $0<\epsilon<\epsilon_{0}$, then\\
$(i)$ (\ref{02}) has a unique nonnegative classical solution $u_{\epsilon}\in C_{loc}^{2,\alpha}(\Omega)$ for some $0<\alpha<1$.\\
$(ii)$ $u_{\epsilon}(x)>0$ for $x\in\Omega$ and
\begin{equation}\label{14011809}
u_{\epsilon}\geq u_{\delta},\;\;\epsilon+u_{\epsilon}\leq \delta+u_{\delta}\;\;\mbox{for}\;\;0<\epsilon\leq\delta\leq\epsilon_{0}.
\end{equation}
\end{lem}
\textbf{Proof.} $(i)$ Similar to \cite{CRT}, in order to establish existence of solution $u_{\epsilon}$, we resort to sub-and supersolutions of (\ref{02}).

Since $(-\Delta)^{s}0=0<\epsilon^{-p}$, then $0$
is a subsolution of (\ref{02}). Next define $w\in C_{loc}^{2,\alpha}(\Omega)$ by
\begin{equation}
\left\{
  \begin{array}{ll}
  (-\Delta)^{s} w=\epsilon^{-p} &\mbox{in}\;\Omega \\
   w=0&\mbox{on}\;\partial\Omega
  \end{array}
\right.
\end{equation}
Then lemmas \ref{mlem14022101}-\ref{mlem14022102} assure the existence and uniqueness of $w$.
By maximum principle, we have $w>0$ in $\Omega$, so, $(-\Delta)^{s}w=\epsilon^{-p}> (\epsilon+w)^{-p}$ in $\Omega$. Thus
$w$ is a supersolution of (\ref{02}) and $[0,w]$ is the desired order pair. Then Lemma \ref{mlem14022110} implies that (\ref{02}) has a solution
$u_{\epsilon}\in C_{loc}^{2,\alpha}(\Omega)$,
with $0<u_{\epsilon}\leq w$ in $\Omega$. The uniqueness follows from $(ii)$.

$(ii)$ Let $0<\epsilon<\delta\leq\epsilon_{0}$. Denote
$u_{\epsilon}$, $u_{\delta}$ such solutions. If $\hat{u}$ is either $u_{\epsilon}-u_{\delta}$ or $(\delta+u_{\delta})-(\epsilon+u_{\epsilon})$,
then $\hat{u}\geq0$ on $\partial \Omega$. (\ref{14011202}) implies, for any constant $C$, $(-\Delta)^{s}C=0$. Let
$$
A=\{x\in\Omega\;|\;\hat{u}<0\}.
$$
Since
$$
(-\Delta)^{s}\hat{u}=(-\Delta)^{s}u_{\epsilon}-(-\Delta)^{s}u_{\delta}\geq0 \;\;\mbox{on}\;A,
$$
and $\hat{u}=0$ on $\partial A$, $A$ must be empty.
Thus $\hat{u}\geq 0$ in $\bar{\Omega}$.

\begin{flushright}
$\Box$
\end{flushright}

 Let $0<\delta_{0}<\delta<\epsilon_{0}$. By Proposition \ref{mpro14030701}, we have
$$
\|u_{\delta}\|_{L^{\infty}(\Omega)}<+\infty,
$$
and so
$$
\|u_{\epsilon}\|_{L^{\infty}(\Omega)}<+\infty.
$$
Then by monotone convergence theorem and Lemma \ref{14031105},
$$
\lim_{\epsilon\rightarrow 0^{+}}u_{\epsilon}=u \;\mbox{exists}\;\mbox{uniformly}\;\mbox{in}\;\bar{\Omega},\;\;\mbox{and}\;\|u\|_{L^{\infty}(\Omega)}<+\infty,
$$
where $u\geq u_{\epsilon} >0$ in $\Omega$ and $u=0$ on $\partial\Omega$.

Let $U_{\epsilon}$ be the extension of $u_{\epsilon}$ in $\mathcal{C}$. Since $u_{\epsilon}\in C_{loc}^{2,\alpha}(\Omega)$, by the boundeness of $u_{\epsilon}$, and $0<p<1$, it is easy to see
\begin{eqnarray*}
&&\int_{\mathcal{C}}y^{1-2s}|\nabla U_{\epsilon}|^{2}dxdy\\
&&=\int_{\Omega\times\{0\}}(\epsilon+u_{\epsilon})^{-p}u_{\epsilon}dx\\
&&\leq\int_{\Omega\times\{0\}}(\epsilon+u_{\epsilon})^{-p}(\epsilon+u_{\epsilon})dx\\
&&=\int_{\Omega\times\{0\}}(\epsilon+u_{\epsilon})^{1-p}dx\\
&&=\int_{\Omega\times\{0\}}(1+u_{\epsilon})^{1-p}dx\\
&&\leq C(\Omega).
\end{eqnarray*}
Then we can extract a subsequence, still denoted by $\{U_{\epsilon}\}$, such that, as $\epsilon\rightarrow 0$,
$$
U_{\epsilon}\rightharpoonup U\;\;\mbox{weakly}\;\mbox{in}\;H^{1}_{0,L}(\mathcal{C}).
$$
By the dominated convergence theorem, for any $\varphi\in H^{1}_{0,L}(\mathcal{C})$, we obtain
\begin{eqnarray*}
&&\int_{\Omega\times\{0\}}u^{-p}\varphi dx=\lim_{\epsilon\rightarrow 0}\int_{\Omega\times\{0\}}(\epsilon+u_{\epsilon})^{-p}\varphi dx\\
&&=\lim_{\epsilon\rightarrow 0}\int_{\mathcal{C}}y^{1-2s}\nabla U_{\epsilon}\nabla \varphi dxdy=\int_{\mathcal{C}}y^{1-2s}\nabla U\nabla \varphi dxdy.
\end{eqnarray*}
Then we deduce that $u$ is a weak solution of (\ref{0310}). By the boundedness of $u$ and Proposition \ref{mpro14030701}, we have
$u\in C_{loc}^{2,\alpha}(\Omega)$ for some $0<\alpha<1$.

Now we prove the uniqueness of $u$. Let $u$ and $v$ be classical solutions of (\ref{0310}). In fact,
$$
u\equiv v.
$$
If not, there exists $x_{0}\in\Omega$ such that
$$
(u-v)(x_{0})=\min_{x\in\Omega}\{(u(x)-v(x))\}<0.
$$
Then
$$
\nabla ((u-v)(x_{0}))=0
$$
However,
$$
y^{1-2s}\frac{\partial (u-v )}{\partial \nu}|_{(x_{0},0)}=\left(u^{-p}(x_{0})-v^{-p}(x_{0})\right)>0,
$$
contradiction.
This completes the proof of Theorem \ref{mthm14011801}.

Authors' Addresses and Emails
\medskip

Yanqin Fang

School of Mathematics

Hunan University

Changsha, 410082 P.R. China, and

Department of Mathematics

Yeshiva University

yanqinfang@126.com
\medskip

\end{document}